\documentclass[a4paper,10pt]{article}

\usepackage[latin1]{inputenc}
\usepackage{amsmath}
\usepackage{amsfonts}
\usepackage{amssymb}
\usepackage[english]{babel}
\usepackage{multirow}
\usepackage[dvips]{graphicx}
\usepackage{color}
\usepackage{indentfirst}
\usepackage{dsfont}
\usepackage{fancyhdr}
\usepackage{setspace}
\usepackage{subfigure}
\usepackage{psfrag}

\makeatletter
\def\thickhrulefill{\leavevmode \leaders \hrule height 1ex \hfill \kern \z@}
\def\@makechapterhead#1{%
  \vspace*{9\p@}%
  {\parindent \z@ \raggedright \reset@font
            \scshape \@chapapp{} \thechapter
        \par\nobreak
        \interlinepenalty\@M
    \Huge \bfseries #1\par\nobreak
    \hrulefill
    \par\nobreak
    \vskip 10\p@
  }}
\def\@makeschapterhead#1{%
  \vspace*{9\p@}%
  {\parindent \z@ \raggedright \reset@font
            \scshape \vphantom{\@chapapp{} \thechapter}
        \par\nobreak
        \interlinepenalty\@M
    \Huge \bfseries #1\par\nobreak
    \hrulefill
    \par\nobreak
    \vskip 10\p@
  }}

\addto{\captionsportuges}{}

\newcommand{\Q}{\mathbb{Q}}
\newcommand{\qed}{\begin{flushright} \vspace{-0.5cm} $\blacksquare$ \end{flushright}}
\newtheorem{tm}{Theorem}[section]
\newtheorem{lemma}[tm]{Lemma}
\newtheorem{defin}[tm]{Definition}
\newtheorem{pp}[tm]{Proposition}

\newtheorem{rem}[tm]{Remark}

\newtheorem{cor}[tm]{Corollary}

\title{Fermat-type equations of signature $(13,13,p)$ via Hilbert cuspforms.}
\author{Luis Dieulefait and Nuno Freitas}

\begin{document}

\maketitle

\begin{abstract} In this paper we prove that equations of the form $x^{13} + y^{13} = Cz^{p}$ have no non-trivial primitive solutions $(a,b,c)$ such that  $13 \nmid c$ if $p > 4992539$ for an infinite family of values for $C$. Our method consists in relating a solution $(a,b,c)$ to the previous equation to a solution $(a,b,c_1)$ of another Diophantine equation with coefficients in $\Q(\sqrt{13})$. We then construct Frey-curves associated with $(a,b,c_1)$  and we prove modularity of them in order to apply the modular approach via Hilbert cusp forms over $\Q(\sqrt{13})$. We also prove a modularity result for elliptic curves over totally real cyclic number fields of interest by itself.
\end{abstract}

\section{Introduction}

Since the proof of Fermat's Last Theorem by Wiles \cite{wiles} the modular approach to Diophantine equations has been popularized and achieved great success in solving equations that previously seemed intractable. In order to attack the generalized Fermat equation $Ax^p + By^q = Cz^r$, where $1/p + 1/q + 1/r < 1$ the initial strategy of Frey, Hellegouarch, Serre, Ribet and Wiles was strengthened and many particular cases, including infinite families, were solved. An important progress is the fact (as a consequence of the work of Darmon-Granville \cite{DG}) that for a fixed triple $(p,q,r)$ there exists only a finite number of solutions such that $(x,y,z)=1$. Another important progress was the work of Ellenberg on the representations attached to $\Q$-curves which allowed to use a special type of elliptic curves over number fields ($\Q$-curves) to attack Diophantine equations over the rationals. In particular, Ellenberg solved the equation $A^4 + B^2 = C^p$ (see \cite{ell}). For a summary of known results on the equation $x^p + y^q = z^r$ see the introduction in \cite{chenBen}.\par
A particularly important subfamily of the generalized Fermat equation are the equations of signature $(r,r,p)$, that is, $Ax^r + By^r = Cz^p$ with $r$ a fixed prime. In this direction there is work for $(3,3,p)$ by Kraus \cite{kraus1}, Bruin \cite{bru}, Chen-Siksek \cite{CS} and Dahmen \cite{Dahm}; for $(5,5,p)$ by Billerey \cite{bil1}, Dieulefait-Billerey \cite{BD} and from the authors \cite{DF}; for $(7,7,p)$ from the second author (currently in preliminary version). In this paper we will go further into this family of equations and we will use a generalized version of the classical modular approach to study equations of the form
\begin{equation}
x^{13} + y^{13} = Cz^p.
\label{treze}
\end{equation}
Let $(a,b,c)$ be a triple of integers such that $a^{13} + b^{13} = Cc^p$. We say that it is a \textit{primitive} solution if $(a,b)=1$ and we will say that it is a \textit{trivial} solution if $abc=0$. Following the terminology introduced by Sophie Germain in her work on the FLT we will divide solutions to (\ref{treze}) into two cases
\begin{defin} A primitive solution $(a,b,c)$ of $x^r + y^r = Cz^p$ is called a first case solution if $r$ do not divide $c$, and a second case solution otherwise.
\end{defin}
The strategy we will use goes as follows: we first relate a possible non-trivial primitive solution $(a,b,c)$ of (\ref{treze}) to a non-trivial primitive solution $(a,b,c_1)$ of another Diophantine equation with coefficients in $\Q(\sqrt{13})$. Secondly we attach to the latter solution a Frey-Hellegouarch-curve over $\Q(\sqrt{13})$ that is not a $\Q$-curve. Then using modularity results from Skinner-Wiles and Kisin we prove modularity of our F-H-curves. With modularity established, from the level lowering results for Hilbert modular forms it will follow that the existence of the solution $(a,b,c_1)$ implies a congruence between two Galois representations. One of the representations is attached to our Frey-curves and the other to a Hilbert newform of a certain level (not depending on $c_1$) and parallel weight $(2,2)$. Finally we will show that, in some cases, this congruence can not hold and so the equation over $\Q(\sqrt{13})$ can not have non-trivial primitive solutions and consequently neither (\ref{treze}). The main result in this paper is the following theorem:

\begin{tm} Let $d=3,5,7$ or $11$ and $\gamma$ be an integer divisible only by primes $l \not\equiv 1$ (mod 13). If $p > 4992539$ is a prime, then:  
\begin{itemize}
\item[(I)] The equation $x^{13} + y^{13} = d\gamma z^p$ has no non-trivial primitive first case solutions.
\item[(II)] The equation $x^{26} + y^{26} = 10\gamma z^p$ has no primitive non-trivial solutions.
\label{teorema}  
\end{itemize}
\label{grande}
\end{tm}
 
In what follows we will first prove part (I) of Theorem \ref{grande} and in the end we will explain the small tweak needed to conclude part (II). Observe that in part (II) replacing $10$ by twice $d$ for $d=3,7,11$ the statement is also true but trivial, because the left-hand side is a sum of two relatively prime squares.\\

We thank John Cremona for providing us a list of elliptic curves that was useful to test our strategy. We also want to thank John Voight for computing a list of Hilbert modular forms that was fundamental to finish this work.

\section{Relating two Diophantine equations}

In this section we will relate a solution of (\ref{treze}) to a solution of a Diophantine equation with coefficients in $\Q(\sqrt{13})$. In order to do that we will need a few properties on the factors of $x^{13} + y^{13}$ over the cyclotomic field $\Q(\zeta_{13})$. Since these properties are not exclusive of degree 13 we will prove them in general. Observe that if $r$ is a prime then $$x^r + y^r = (x+y)\phi_r(x,y)$$ where $$\phi_r (x,y)= \prod_{i=0}^{r-1}{(-1)^{i}x^{r-1-i}y^i}.$$
Let $\zeta = \zeta_r$ be an $r$-root of unity and consider the decomposition over the cyclotomic field $\Q(\zeta)$
\begin{equation}
\label{decomp}
\phi_r (x,y) = \prod_{i=1}^{r-1}{(x + \zeta^i y)}.
\end{equation}

\begin{pp} Let $\mathfrak{P}_r$ be the prime in $\Q(\zeta)$ above the rational prime $r$ and suppose that $(a,b)=1$. Then, any two distinct factors $a + \zeta^i b$ and $a + \zeta^j b$ in the factorization of $\phi_r (a,b)$ are coprime outside $\mathfrak{P}_r$. Furthermore, if $r \mid a+b$ then $\nu_{\mathfrak{P}_r}(a + \zeta^i b) = 1$ for all $i$.
\label{factores}
\end{pp}
\textbf{Proof:} Suppose that $(a,b)=1$. Let $\mathfrak{P}$ be a prime in $\Q(\zeta)$ above $p \in \Q$ and a common prime factor  of $a + \zeta^i b$ and $a + \zeta^j b$, with $i >j$. Observe that $(a + \zeta^i b) - (a + \zeta^j b) = b\zeta^j(1-\zeta^{i-j})$ and since $\mathfrak{P}$ must divide the difference it can not divide $b$ because in this case it would also divide $a$, and since $a,b$ are integers $p$ would divide both. So $\mathfrak{P}$ must be a factor of $\zeta^i(1-\zeta^{i-j})$ but $\zeta^{i}$ is a unit for all $i$ then $\mathfrak{P}$ divides $1-\zeta^{i-j}$, that is $\mathfrak{P} = \mathfrak{P}_{r}$. Now for the last statement in the proposition, suppose that $r \mid a+b$. Then,
$$ a + \zeta^i b = a + b - b + \zeta^i b = (a+b) + (\zeta^i - 1)b,$$ 
and since $\nu_{\mathfrak{P}_r}(\zeta^i - 1) = 1$ we have $\nu_{\mathfrak{P}_r}(a + \zeta^i b) = \mbox{min}\{r-1,1\} = 1$\qed

\begin{cor} If $(a,b)=1$, then $a+b$ and $\phi_r (a,b)$ are coprime outside $r$. Furthermore, if $r \mid a+b$ then $\nu_{r}(\phi_r (a,b))=1$. 
\label{trezz}
\end{cor}
\textbf{Proof:} Let $p$ be a prime dividing $a+b$ and $\phi_r (a,b)$ and denote by $\mathfrak{P}$ a prime in $\Q(\zeta)$ above $p$. $\mathfrak{P}$ must divide at least one of the factors $a + \zeta^i b$. Since $a,b$ are integers $\mathfrak{P}$ can not divide $b$ then it follows from $$a+b = a + \zeta^i b - \zeta^i b + b = (a + \zeta^i b) + (1-\zeta^i)b$$ that $\mathfrak{P} = \mathfrak{P}_r$. Moreover, if $r \mid a+b$ it follows from the proposition that $\nu_{\mathfrak{P}_r}(a + \zeta^i b) = 1$ for all $i$ then $\nu_{\mathfrak{P}_r}(\phi_r(a,b)) = r-1$ thus $\nu_{r}(\phi_r (a,b)) = 1$.\qed

\begin{pp} Let $(a,b)=1$ and $l \not\equiv 1$ (mod $r$) be a prime dividing $a^{r} + b^{r}$. Then $l \mid a+b$.
\label{trezz2}
\end{pp}
\textbf{Proof:} Since $l$ divides $a^{r} + b^{r}$, $l \nmid ab$. Let $b_0$ be the inverse of $-b$ modulo $l$. We have $a^{r} \equiv (-b)^{r}$ (mod $l$), hence $(ab_0)^{r} \equiv 1$ (mod $l$). Thus the multiplicative order of $ab_0$ in $\mathbb{F}_l$ is 1 or $r$. From the congruence $ab_0 \equiv 1$ (mod $l$) it follows $a + b \equiv 0$ (mod $l$). If $l \nmid a+b$ then the order of $ab_0$ is $r$ and $l \equiv 1$ (mod $r$).\qed

From now on we particularize to $r=13$ and we denote $\phi_{13}$ only by $\phi$. We have $x^{13} + y^{13} = (x + y)\phi(x,y)$, where
\begin{eqnarray*}
\phi(x,y) & = &  x^{12} - x^{11}y + x^{10}y^2 - x^9y^3 + x^8y^4 - x^7y^5 + x^6y^6 \\
          & - &  x^5y^7 + x^4y^8 - x^3y^9 + x^2y^{10} - xy^{11} + y^{12}. \
\end{eqnarray*}
Suppose that there exists a non-trivial primitive solution $(a,b,c)$ to (\ref{treze}) with $C = d \gamma$, $d$ and $\gamma$ as in theorem \ref{grande}. Then it follows from corollary \ref{trezz} and proposition \ref{trezz2} that there exists a non-trivial primitive solution $(a,b,c_0)$ to
\begin{equation}
\phi(a,b) = c_0^p,
\label{casoA}  
\end{equation}
with $d \mid a+b$ and $13 \nmid a+b$ or to
\begin{equation}
\phi(a,b) = 13c_0^p
\label{casoB} 
\end{equation}
with $d \mid a+b$ and $13 \mid a+b$, where in both cases $c_0$ is only divisible by primes congruent to 1 modulo 13. Consider the factorization of $\phi$ in $\Q(\zeta)$
$$\phi(x,y) = \prod_{i=1}^{12}{(x+\zeta^{i}y)}.$$
We have $\phi = \phi_1 \phi_2$ where $\phi_i$ are both of degree 6 with coefficients in $\Q(\sqrt{13})$, given by  
\begin{eqnarray*}
\label{phifactors}
\phi_1(x,y) & = & (x+\zeta y)(x+\zeta^{12} y)(x+\zeta^{4} y)(x+\zeta^{9} y)(x+\zeta^{3} y)(x+\zeta^{10} y),\\
\phi_2(x,y) & = & (x+\zeta^2 y)(x+\zeta^{5} y)(x+\zeta^{6} y)(x+\zeta^{7} y)(x+\zeta^{8} y)(x+\zeta^{11} y).\
\end{eqnarray*}

The proof of the following corollary of proposition \ref{factores} is immediate.
\begin{cor} Let $\mathfrak{P}_{13}$ be the prime of $\Q(\sqrt{13})$ above 13. If $(a,b) = 1$, then $\phi_1(a,b)$ and $\phi_2(a,b)$ are coprime outside $\mathfrak{P}_{13}$. Moreover, $\nu_{\mathfrak{P}_{13}}(\phi_i(a,b)) = 1$ or 0 for both $i$ if $13 \mid a+b$ or $13 \nmid a+b$, respectively. 
\label{factores2}
\end{cor}

Corollary \ref{factores2} and the existence of a solution $(a,b,c_0)$ to (\ref{casoA}) or (\ref{casoB}) implies that for some unit $\mu$ there exists a solution $(a,b,c_1)$ (with $c_1$ an integer  in $\Q(\sqrt{13})$) to the equation 
\begin{equation}
\label{caso1}
\phi_1(a,b) = \mu c_1^p, 
\end{equation}
with $d \mid a+b$ and $13 \nmid a+b$ or to
\begin{equation}
\phi_1(a,b) =\mu \sqrt{13}c_1^p, 
\label{caso13}
\end{equation}
with $d \mid a+b$ and $13 \mid a+b$, respectively.\\
 
Observe that proposition \ref{trezz2} and the form of equation (\ref{treze}) implies that $13 \mid c$  is equivalent to $13 \mid a+b$. Moreover, proposition \ref{trezz2} also guarantees that when passing from the equation in (I) of Theorem \ref{grande} to equations (\ref{caso1}) or (\ref{caso13}) the prime factors of $\gamma$ can be supposed to divide $a+b$. Since this information will not be necessary to the proof of both parts of theorem \ref{grande} we can assume that $\gamma = 1$. It will also be clear from the proof that the unit $\mu$ can be supposed to be 1.  

\section{The Frey-Hellegouarch curves}

Let $\sigma$ be the generator of $G=\mbox{Gal}(\Q(\zeta)/\Q)$ and $K$ the subfield (of degree 6) fixed by $\sigma^6$. Consider the factorization $\phi_1 = f_1 f_2 f_3$ where 
$$\begin{cases}
 f_1(x,y) = (x+\zeta y)(x+\zeta^{12} y)=x^2 + (\zeta + \zeta^{12})xy + y^2 \\
 f_2(x,y) = (x+\zeta^{4} y)(x+\zeta^{9} y)=x^2 + (\zeta^4 + \zeta^{9})xy + y^2\\
 f_3(x,y) = (x+\zeta^{3} y)(x+\zeta^{10} y)=x^2 + (\zeta^3 + \zeta^{10})xy + y^2\\
\end{cases}$$
are the degree two factors of $\phi_1$ with coefficients in $K$. Now we are interested in finding a triple $(\alpha, \beta, \gamma)$ such that 
$$\alpha f_1 + \beta f_2 + \gamma f3 = 0.$$ 
Solving a linear system in the coefficients of the $f_i$ we find that one of its infinite solutions in $\mathcal{O}_K^3$ is given by
$$\begin{cases}
 \alpha = -\zeta^{10} + \zeta^9 + \zeta^4 - \zeta^3 \\
 \beta  = \zeta^{12} - \zeta^9 - \zeta^4 + \zeta \\
 \gamma =  -\zeta^{12} + \zeta^{10} + \zeta^3 - \zeta\\
\end{cases}$$
and verifies $\nu_{\mathfrak{P}_{13}}(\alpha) = \nu_{\mathfrak{P}_{13}}(\beta) = \nu_{\mathfrak{P}_{13}}(\gamma) = 1$.\par 
Suppose now that $(a,b,c) \in \mathbb{Z} \times \mathbb{Z} \times \mathcal{O}_{\Q(\sqrt{13})}$ is a non-trivial primitive solution to equation (\ref{caso1}) or (\ref{caso13}) and let $A(a,b) = \alpha f_1(a,b)$, $B(a,b) = \beta f_2(a,b)$ and $C(a,b) = \gamma f_3(a,b)$. Since  $$A + B + C = 0$$ then we can consider the Frey-curve over $K$ with the classic form
$$E(a,b): y^2 = x(x-A(a,b))(x+B(a,b)).$$
In the rest of this section we will denote $E(a,b)$ only by $E$ every time it causes no ambiguity. Let $\mathfrak{P}_{2}$ denote the prime of $K$ above 2 (is inert). To the curves $E(a,b)$ are associated the following quantities:
\begin{eqnarray*}
\Delta(E) & = & 2^{4}(ABC)^2 = 2^4(\alpha \beta \gamma)^2\phi_1(a,b)^2  \\ 
c_4(E) & = & 2^{4}(A^2 + AB + B^2) = 2^{4}(AB + BC + AC) \\
c_6(E) & = & -2^{5}(C + 2B)(A + 2B)(2A + B) \\
j(E) & = & 2^{8}\frac{(A^2 + AB + B^2)^3}{(ABC)^2} \
\end{eqnarray*}
and since $(\alpha\beta\gamma) = \mathfrak{P}_{13}^3$ the discriminant of $E$ takes the following values  
$$ \Delta(E) = \begin{cases} 
\mathfrak{P}_{2}^{4}\mathfrak{P}_{13}^{6}c^{2p} \mbox{ if } 13  \nmid a+b \\ 
\mathfrak{P}_{2}^{4}\mathfrak{P}_{13}^{12}c^{2p} \mbox{ if } 13 \mid a+b
\end{cases}$$

\begin{pp} Let $\mathfrak{P}$ be a prime of $K$ distinct from $\mathfrak{P}_2$ and $\mathfrak{P}_{13}$. The curves $E(a,b)$ have good or multiplicative reduction at $\mathfrak{P}$. Moreover, the curves have good ($\nu_{\mathfrak{P}_{13}}(N_E) = 0$) or bad additive ($\nu_{\mathfrak{P}_{13}}(N_E) = 2$) reduction at $\mathfrak{P}_{13}$ if $13 \mid a+b$ or $13 \nmid a+b$, respectively. In particular, $E$ has multiplicative reduction at primes dividing $c$. 
\label{condE}
\end{pp}

\textbf{Proof:} To the results used in this proof we followed \cite{pap}. Let $\mathfrak{P}$ be has in the hypothesis and observe that 
$\upsilon_{\mathfrak{P}}(\Delta(E)) = 2p\upsilon_{\mathfrak{P}}(c)$. Then if $\mathfrak{P} \nmid c$ we have $\upsilon_{\mathfrak{P}}(\Delta) = 0$ and the curve has good reduction. It follows from proposition \ref{factores} that $A(a,b)$, $B(a,b)$ and $C(a,b)$ are pairwise coprime outside $\mathfrak{P}_{13}$ and recall that the three are divisible by $\mathfrak{P}_{13}$. If $\mathfrak{P} \mid c$ then $\mathfrak{P}$ must divide only one among $A,B$ or $C$. From the form of $c_4$ it can be seen that $\upsilon_{\mathfrak{P}}(c_4)=0$. Also, $\upsilon_{\mathfrak{P}}(\Delta) > 0$ thus $E$ has multiplicative reduction at $\mathfrak{P}$. Moreover, we see from $\Delta(E)$ that if $13 \mid a+b$ the equation is not minimal and $E$ has good reduction at $\mathfrak{P}_{13}$. On the other hand if $13 \nmid a+b$ then $\nu_{\mathfrak{P}_{13}}(\Delta) = 6$ and $\nu_{\mathfrak{P}_{13}}(c_4) > 0$ hence bad additive reduction at  $\mathfrak{P}_{13}$.\qed
 
\begin{pp} The short Weierstrass model of the curves $E(a,b)$ is defined over $\Q(\sqrt{13})$.
\end{pp}
\textbf{Proof:} First observe that $\sigma^{4}$ (mod $\sigma^6$) generates $\mbox{Gal}(K/\Q(\sqrt{13}))$. Since the curves $E$ are defined over $K$ they are invariant under $\sigma^6$ and in particular $j(E)$ is invariant  by $\sigma^6$ by definition. We also have that 
$$ \sigma^{4}(A) = B , \hspace{1cm} \sigma^{4}(B) = C , \hspace{1cm}  \sigma^{4}(C) = A,$$  
and from 
$$j(E) =  2^{8}\frac{(AB + BC + CA)^3}{(ABC)^2}$$ 
it is clear that $j$ is also invariant under $\sigma^{4}$. Then the $j$-invariant that \textit{a priori} belonged to $K$ of degree $6$, in reality is in $\Q(\sqrt{13})$. Now we write $E(a,b)$ in the short Weierstrass form to get a model
$$\begin{cases}
   E_0 : y^2 = x^3 + a_4 x + a_6, \mbox{ where } \\
   a_4 = -432(AB + BC + CA) \\
   a_6 = -1728 (2A^3 + 3A^2B - 3AB^2 - 2B^3) \\
\end{cases}$$
Since $a_4$ is clearly invariant under $\sigma^{4}$ and
\begin{eqnarray*}
  a_6  & = & -1728 (2A^3 + 3A^2B - 3AB^2 - 2B^3) = \\
       & = &  -1728 (2(-B-C)^3 + 3(-B-C)^2B - 3(-B-C)B^2 - 2B^3) = \\
       & = & -1728 (2B^3 + 3B^2C - 3BC^2 - 2C^3) = \sigma^{4}(a_6)  \
\end{eqnarray*}
we conclude that the short Weierstrass model is already defined over $\Q(\sqrt{13})$.\qed

Writing $E$ in the short Weierstrass form we get an elliptic curve $E_{0}$ defined over $\Q(\sqrt{13})$ given by 
$$E_{0} : y^2 = x^3 + a_4(a,b)x + a_6(a,b),$$
\begin{eqnarray*}
 a_4(a,b) & = &  (216w - 2808)a^4 + (-1728w + 5616)a^3b \\
	  &   &  + (1728w - 11232)a^2 b^2 + (-1728w + 5616)ab^3 \\
	  &   &  + (216w - 2808)b^4 ,  \\
 a_6(a,b) & = &  (-8640w + 44928)a^6 + (49248w - 235872)a^5b \\
          &   &  + (-129600w + 471744)a^4b^2 + (152928w - 662688)a^3b^3  + \\ 
          &   &  + (-129600w + 471744)a^2b^4 + (49248w - 235872)ab^5  + \\
	  &   &  + (-8640w + 44928)b^6 + (50193w + 182520)b^6 ,\ 
\end{eqnarray*}
where $w^2 = 13$. Writing a curve in short Weierstrass form changes the values of $\Delta$, $c_4$ and $c_6$ according to $\Delta(E_0) = 6^{12}\Delta(E)$, $c_4(E_{0}) = 6^4 c_4(E)$ and $c_6(E_{0}) = 6^6 c_6(E)$. Observe that 2 is inert in $\Q(\sqrt{13})$ (and in $K$) and denote by 2 and $w$ the ideals in $\Q(\sqrt{13})$ above 2 and 13, respectively. 

\begin{pp} The possible values for the conductors of $E_{0}(a,b)$ are 
$$N_{E_0} = 2^{s}w^{2}\mbox{rad}(c),$$
where $s=3,4$ and $\mbox{rad}(c)$ is the product of the prime factors of $c$. Moreover, if $2 \mid a+b$ then $s=3$ if  $4 \mid a+b$ and $s=4$ if $4 \nmid a+b$. 
\label{condEgamma} 
\end{pp}
\textbf{Proof:} As in proposition \ref{condE} we followed the results in \cite{pap} to compute the conductor. Since the primes dividing $6$ do not ramify in $K/\Q(\sqrt{13})$ and do not divide $c$ the conductor of $E_{0}$ and $E$ is same at these primes. \par 
Since $(w) = \mathfrak{P}_{13}^3$ in $K$ we see from Proposition \ref{condE} that $\nu_{w}(\Delta(E_0)) = 4$ or $2$ if $13 \mid a+b$ or $13 \nmid a+b$, respectively. Also, $\nu_{w}(c_4(E_0)) > 0$ and since we are in characteristic $\geq 5$ this implies that the equation is minimal and has bad additive reduction with $\nu_w (N_{E_{0}}) = 2$.\par
It easily can be seen that $\nu_{2}(\Delta(E_0)) = 16$, $\nu_{2}(c_6(E_0)) = 11$ and $\nu_{2}(c_4(E_0)) \geq 8$. Table IV in \cite{pap} tell us that the equation is not minimal and after a change of variables we have $\nu_{2}(\Delta(E_0)) = 4$, $\nu_{2}(c_6(E_0)) = 5$ and $\nu_{2}(c_4(E_0)) \geq 4$. We check that in the columns where $\nu_{2}(c_6(E_{0})) = 5$ and $\nu_{2}(\Delta(E_{0})) = 4$ we have $\nu_{2}(N_{E_{0}}) = 2, 3, 4$.
The conductors at 2 of the curves $E(1,-1)$ and $E(1,1)$ are $2^3$ and $2^4$, respectively. The case $s=2$ never happens. This is a direct consequence of Proposition 2 in \cite{pap} by taking $r=w/2 + 1/2$ if $(a,b) \equiv (1,0),(0,1) \mbox{ (mod } 2)$ or $r=0$ if $(a,b) \equiv (1,1)$ (mod $2$). Since the same proposition 2 guarantees that all the possible conductors at 2 will occur for the pairs $(a,b)$ modulo 4 we use SAGE to compute the conductor for all these pairs and easily verify the statement by inspection.\qed 

From now on we will write $E$ to denote $E_0$.

\section{The Galois representations of $E(a,b)$}

Let $\rho_{E,p} : G_{\Q(\sqrt{13})} \rightarrow GL_{2}(\Q_p)$ be the $p$-adic representation associated with $E$ and $\bar{\rho}_{E,p}$ its reduction modulo $p$. In this section we will prove that $\rho_{E,p}$ is modular and that $\bar{\rho}_{E,p}$ is irreducible for $p>97$. These results together allow us to apply the lowering the level theorems for Hilbert modular forms.

\subsection{Modularity}

\begin{tm} Let $F$ be a totally real cyclic number field and $E$ and elliptic curve defined over $F$. Suppose that 3 splits in $F$ and $E$ has good reduction at the primes above 3. Then $E$ is modular.
\label{modularidade}
\end{tm}

Let $N_E$ denote the conductor of $E$ and put $\rho = \rho_{E,3}$. Since the representation $\bar{\rho} = \bar{\rho}_{E,3} : \mbox{Gal}(\bar{\Q}/\Q(\sqrt{13})) \rightarrow \mbox{GL}_2(\mathbb{F}_3)$ is odd it is absolutely irreducible if and only if it is irreducible. The following key lemma is a known result and a consequence of the work of Langlands-Tunnell. For references where it is used see \cite{mano}, \cite{taylor}, \cite{el2} and \cite{JMano}. 
   
\begin{lemma} If $\bar{\rho}$ is irreducible then it is modular arising from an Hilbert newform over $F$ of parallel weight $(2,2)$. Moreover, if $\rho$ is ordinary we can suppose that $f$ is ordinary at 3.
\label{resmod}
\end{lemma}
\textbf{Proof:} Let $t$ be a prime in $F$ above 3.\par 
For the second statement of the theorem, by hypothesis we know that 
$$ {\bar{\rho}}|{I_{t}} = \begin{pmatrix} \bar{\chi} & * \\ 0 & 1 \end{pmatrix},$$
where $\chi$ is the 3-adic cyclotomic character. It is know from the work of Jarvis, Rajaei and Fujiwara that since $\bar{\rho}_{f,3}|D_{t} \equiv \bar{\rho}|D_{t}$ and $\rho|D_{t}$ is a Barsotti-Tate representation for all $t$ (because $E$ has good reduction at 3) we can choose $f$ to have parallel weight $(2,2)$ and level coprime with 3, hence $\rho_{f,3}|D_{t}$ is also Barsotti-Tate. Since $\bar{\rho}_{f,3}|D_{t}$ is ordinary we can apply a result of Breuil (see \cite{Bre}, chapter 9) to conclude that $\rho_{f,3}$ is ordinary.\qed

We will now prove the theorem. A similar argument but over $\Q$ was given by the first author in \cite{Die}.\par

\textbf{Proof} (of theorem \ref{modularidade}): We divide the proof into three cases:

\begin{itemize}
\item[(1)] Suppose that $\bar{\rho}$ and ${\bar{\rho}}|G_{\Q(\sqrt{-3})}$ are both abs. irreducible. Here we apply corollary 2.1.3 in \cite{kis2}. Condition (1) holds because $E$ has good reduction at the primes above 3 and 3 splits in $F$. Lemma \ref{resmod} guarantees condition (2) and (3) is obvious. Then $\rho$ is modular.
\item[(2)] Suppose that $\bar{\rho}$ is abs. irreducible and ${\bar{\rho}}|G_{F(\sqrt{-3})}$ abs. reducible. This means that the image of $\mathbb{P}(\bar{\rho})$ is Dihedral. Namely, that the image of $\bar{\rho}$ is contained in the normalizer $N$ of a Cartan subgroup $C$ of $GL_2(\bar{\mathbb{F}}_3)$ but not contained in $C$. Moreover, the restriction to $\Q(\sqrt{-3})$ of our representation has its image inside $C$. Thus, the composition of $\bar{\rho}$ with the quotient $N/C$,
\begin{equation}
\label{comp}
\mbox{Gal}(\bar{\Q}/F) \rightarrow N \rightarrow N/C, 
\end{equation}
gives the quadratic character of $F(\sqrt{-3}) / F$ which ramifies at 3 because 3 is unramified in $F$.\par

Let $t$ be a prime in $F$ above 3. Since $E$ has good reduction at 3, the restriction of the residual representation $\bar{\rho}$ to the inertia subgroup $I_t$ has only two possibilities
$$ {\bar{\rho}}|{I_{t}} = \begin{pmatrix} \bar{\chi} & * \\ 0 & 1 \end{pmatrix} \quad \mbox{ or } \quad 
\begin{pmatrix} \psi_2 & 0 \\ 0 & \psi_2^3 \end{pmatrix},$$
where $\chi$ is the 3-adic cyclotomic character and $\psi_2$ is the fundamental character of level 2. \par
If we suppose that $\bar{\rho}|I_t$ acts through level 2 fundamental characters, the image of $I_t$ by $\mathbb{P}(\bar{\rho})$ gives a cyclic group of order $4 > 2$, thus it has to be contained in $\mathbb{P}(C)$ (if it not contained in $\mathbb{P}(C)$ and has order 4 it must be  isomorphic to $C_2 \times C_2$). But this implies that the quadratic character defined by composition (\ref{comp}) should be unramified at 3, contradicting the fact that this character corresponds to $F(\sqrt{-3})$. Thus we can assume that we are in the first case, that is, $\bar{\rho}|I_{t}$ is reducible. Since $\rho$ is crystalline with Hodge-Tate weights 0 and 1 and 3 splits in $F$ we apply a result of Breuil (see \cite{Bre}, chapter 9) to conclude that $\rho|D_{t}$ is reducible hence ordinary. Indeed, we can apply this result of classification of
crystalline representations at $3$ because the highest Hodge-Tate weight is $w = 1$ and so  $3 > w + 1$. Since $t$ is arbitrary the previous holds for all primes $t$ above 3 hence $\rho$ is ordinary. Therefore, we can suppose that the form $f$ given by Lemma \ref{resmod} is ordinary thus all conditions are satisfied to apply  Theorem 5.1 in \cite{SW2} to conclude that $\rho$ is modular. 
\item[(3)] Suppose that $\bar{\rho}$ is abs. reducible. We exclude again the case of the fundamental characters of level 2, but this time this is automatic 
because of reducibility. Then
$$ \bar{\rho}^{ss} = \epsilon \oplus \epsilon^{-1}\bar{\chi}, $$
where $\epsilon$ ramifies only at the primes of $N_E$. Also, since the representation is odd and $F$ is real, reducibility must take place over $\mathbb{F}_3$.   Again by the result of Breuil on crystalline representations we conclude that $\rho|D_{t}$ is reducible. Observe now that $\epsilon$ must be quadratic because $\mathbb{F}_3^{*}$ has two elements then $\epsilon/(\epsilon^{-1}\bar{\chi}) = \bar{\chi}$ and the extension $F(\sqrt{-3})/\Q$ is abelian because $F$ is cyclic and disjoint from $\Q(\sqrt{-3})$. This establishes all the conditions of theorem A in \cite{SW1} thus $\rho$ is  modular. \qed
\end{itemize}

\subsection{Irreducibility}

\begin{pp}
Let $p > 97$ be a prime. The representation $\bar{\rho}_{E,p}$ is absolutely irreducible.
\end{pp}
\textbf{Proof:} Since $\bar{\rho}_{E,p}$ is odd and $\Q(\sqrt{13})$ is totally real it is known that $\bar{\rho}_{E,p}$ is  absolutely irreducible if and only if it is irreducible then we only need to rule out the case where $\bar{\rho}_{E,p}$ is reducible and has the form 
\begin{equation}
\label{red}
\bar{\rho}_{E,p} = \begin{pmatrix} \epsilon^{-1}\chi_{p} & * \\ 0 & \epsilon \end{pmatrix},
\end{equation} 
where $\chi_p$ is the mod $p$ cyclotomic character and $\epsilon$ is a character of $G_{\Q(\sqrt{13})}$ with values in $\mathbb{F}_{p}$. Since the image of inertia at semistable primes is of the form $\begin{pmatrix} 1 & * \\ 0 & 1 \end{pmatrix}$ the conductor of $\epsilon$ only contains additive primes. By the work of Carayol the conductor at bad additive primes of $\bar{\rho}_{E,p}$ is the same as that of $\rho_{E,p}$. Since the conductors of $\epsilon$ and $\epsilon^{-1}$ are equal it follows from proposition \ref{condEgamma} that the $\mbox{cond}(\epsilon)=\mathfrak{P}_{2}\mathfrak{P}_{13}$ or $\mathfrak{P}_{2}^2\mathfrak{P}_{13}$. The characters of $G_{\Q(\sqrt{13})}$ with conductor dividing $\mathfrak{P}_{2}^2 \mathfrak{P}_{13}$ are in correspondence with the characters of the finite group 
$$H = (\mathcal{O}_{\Q(\sqrt{13})} / \mathfrak{P}_{2}^2 \mathfrak{P}_{13})^{*} \cong \mathbb{Z}/12\mathbb{Z} \times \mathbb{Z}/6\mathbb{Z} \times \mathbb{Z}/2\mathbb{Z}.$$
The group of characters of $H$ is dual of $H$ then all the characters have order dividing 12. In particular $\epsilon$ is a root of the polynomial $q_1 := x^{12}-1$ (mod $p$) for any $p$. Let $\mathfrak{P}_3$ be a prime above 3. By taking traces on equality (\ref{red}) we get
$$ a_{\mathfrak{P}_3} \equiv \epsilon(\mbox{Frob}_{\mathfrak{P}_3}) + 3\epsilon^{-1}(\mbox{Frob}_{\mathfrak{P}_3}) \hspace{1cm} (\mbox{mod p}),$$
which implies that $\epsilon(\mbox{Frob}_{\mathfrak{P}_3})$ satisfies for any $p$ the polynomial $q_2 := x^2 - a_{\mathfrak{P}_3}x + 3$ (mod $p$). Let $\zeta=\zeta_{12}$, then the resultant of $q_1$ and $q_2$ is given by
\begin{eqnarray*}
 \mbox{res}(q_1,q_2) & = & \prod_{i=1}^{12}{(\frac{a_{\mathfrak{P}_3} + \sqrt{a_{\mathfrak{P}_{3}}^2} - 12}{2} - \zeta^{i})(\frac{a_{\mathfrak{P}_3} - 
  \sqrt{a_{\mathfrak{P}_{3}}^2} - 12}{2} - \zeta^{i})} \\
		    & = & \prod_{i=1}^{12}({\zeta^{2i} - a_{\mathfrak{P}_{3}}\zeta^i + 3})
\end{eqnarray*}
Since $|a_{\mathfrak{P}_{3}}| \leq 3$ we compute all the possibilities to the product above to find that the greater prime divisor appearing as a possible factor is 97. Therefore, since $\epsilon(\mbox{Frob}_{\mathfrak{P}_3})$ is a common root of the $q_i$ (mod $p$) then $\mbox{res}(q_1,q_2) \equiv 0$ (mod $p$), which is impossible if $p > 97$. Thus $\bar{\rho}_{E,p}$ is absolutely irreducible if $p>97$.\qed

With the theorems above we are now able to lower the level. See Jarvis-Meekin \cite{JM} for an application of the level lowering results for Hilbert modular forms in \cite{Raj}, \cite{Fuj} and \cite{Jarv}. In the present case we apply these results along the same lines. Denote by $S_2 (N)$ the set of Hilbert modular cusp forms of parallel weight $(2,2)$ and level $N$. It follows from the modularity that there exists a newform $f_0$ in $S_2 (\mathfrak{P}_{2}^{i}\mathfrak{P}_{13}^{2}\mbox{rad}(c))$ where $i=3$ or $4$ such that $\rho_{E,p}$ is isomorphic to the $p$-adic representation associated with $f_0$, which we denote by $\rho_{f_0,p}$. Since the semistable primes of $E$, i.e. those dividing $c$, appear to a $p$-th power in the discriminant $\Delta(E)$ we know by an argument of Hellegouarch that the representation $\bar{\rho}_{E,p}$ will not ramify at these primes. Furthermore, when reducing to the residual representation the conductor at the bad additive primes can not decrease hence $\bar{\rho}_{E,p}$ has conductor equal to that of $\rho_{E,p}$ without the factor $\mbox{rad}(c)$, that is $\mathfrak{P}_{2}^{i}\mathfrak{P}_{13}^{2}$ with $i=3$ or $4$. Since $\rho_{E,p}$ is modular then $\bar{\rho}_{E,p}$, when irreducible, is modular and by the results on level lowering for Hilbert modular forms  we know that there exists a newform $f$ in $S_2 (\mathfrak{P}_{2}^{i}\mathfrak{P}_{13}^{2})$ such that its associated mod $p$ Galois representation satisfies 
\begin{equation}
\label{congruencia}
\rho_{E,p} \equiv \rho_{f_0,p} \equiv \rho_{f,p} \mbox{ (mod } \mathfrak{P}). 
\end{equation}

\section{Eliminating Newforms}

In this section we will find a contradiction to congruence (\ref{congruencia}). This shows that the Frey-curves associated with primitive non-trivial first case solutions $(a,b,c)$ to equation (\ref{caso1}) or (\ref{caso13}) can not exist and ends the proof of part (I) in Theorem \ref{grande}. To find the desired contradiction we use the trace values $a_{L}(\rho_{E,p})$ and $a_{L}(\rho_{f,p})$ for some primes $L$ of $\Q(\sqrt{13})$ and the Hilbert modular newforms $f$ in the spaces predicted in the previous section. Let $w \in \Q(\sqrt{13})$ be such that $w^2 = 13$ and consider the following prime ideals in $\Q(\sqrt{13})$:
$$\begin{cases}
L_{2} = \langle 2 \rangle, \quad L_{13} = \langle w \rangle \\
L_{3}^{0} =  \langle \frac{1}{2}(w + 1) \rangle, \quad  L_{3}^{1} = \langle \frac{1}{2}(-w + 1)\rangle, \\ 
L_{17}^{0} = \langle \frac{1}{2}(w + 9)\rangle, \quad  L_{17}^{1} = \langle \frac{1}{2}(-w + 9)\rangle, \\
L_{23}^{0} = \langle \frac{1}{2}(-3w - 5)\rangle, \quad L_{23}^{1} = \langle \frac{1}{2}(-3w + 5)\rangle, \\
L_{29}^{0} = \langle \frac{1}{2}(3w + 1) \rangle, \quad L_{29}^{1} = \langle \frac{1}{2}(3w - 1)\rangle, \\
L_{5} = \langle 5 \rangle, \quad L_{7} = \langle 7 \rangle, \quad L_{11} = \langle 11 \rangle. \
\end{cases}$$

On one hand, to obtain the values of $a_L(\rho_{f,p})$, with the aid of John Voight we used algorithms to compute Hilbert modular forms implemented in MAGMA \cite{magma} (an expository account can be found in \cite{DV}). John Voight gave us two lists corresponding to all forms with integer coefficients such that $a_{L_2} = a_{L_{13}} = 0$ and of levels $\mathfrak{P}_{2}^{s}\mathfrak{P}_{13}^{2}$ for $s=3,4$. With MAGMA we have done the same to all dividing levels and by putting together both informations we obtained all newforms in the spaces $S_2(\mathfrak{P}_{2}^{s}\mathfrak{P}_{13}^{2})$ for $s=3,4$ such that $\Q_f =\Q$. A list of coefficients corresponding to the newforms obtained this way can be found in the appendix A. Moreover, a consequence of the method used is that any newform in the two previous spaces with $\Q_f$ strictly containing $\Q$ must have a Fourier coefficient outside $\Q$ at the prime $L_{3}^{0}$ above 3. John Voight also computed the factorization of the characteristic polynomial of the Hecke operator $T_{L_3^{0}}$ in both spaces (see appendix B).\\

On the other hand, for every prime $L$ in $\Q(\sqrt{13})$ of good reduction for $E$, such that $L$ is above a rational prime $l \leq 29$ and $l \not= 19$, we use SAGE to go through all the possible residual elliptic curves for all pairs $(a,b) \in \mathbb{F}_l \times \mathbb{F}_l$  and compute all the possible values for $a_{L}(\rho_{E,p}) = a_L(E)$: 
$$\begin{cases} 
 a_{L_3^{0}} \in \{-3,-1\}, \\
 a_{L_{3}^{1}} \in \{-3,-1,1\}, \\ 
 a_{L_{5}} \in \{-6, -2 ,2\}, \\ 
 a_{L_{7}} \in \{11, -11, -1, -5\}, \\
 a_{L_{11}} \in \{ -15, 3, 5, -7, 9, -1, 15 \}, \\
 a_{L_{17}^{0}} \in \{1, 3, 5, 7, -3, -1\}, \\
 a_{L_{17}^{1}} \in \{3, 5, 7, -7, -5, -3 \}, \\ 
 a_{L_{23}^{0}} \in \{1, 3, 5, 7, -9, -7, -5, -3\}, \\
 a_{L_{23}^{1}} \in \{1, 3, 7, -9, -3, -1\}, \\
 a_{L_{29}^{0}} \in \{1, 3, 5, -9, -7, -5, -3, -1\}, \\
 a_{L_{29}^{1}} \in \{1, 3, 5, 9, -9, -7, -5, -3, -1\} 
\end{cases}$$
Before proceeding to eliminate the newforms we divide them into two sets: 
\begin{itemize}
 \item S1: The newforms in $S_2 (\mathfrak{P}_{2}^{i}\mathfrak{P}_{13}^{2})$ for $i=3,4$ such that $\Q_f = \Q$. 
 \item S2: The newforms in the same levels with $\Q_f$ strictly containing $\Q$.
\end{itemize}
Note that equations (\ref{caso1}) and (\ref{caso13}) have trivial solutions $(1,1,1)$, $\pm (0,1,1)$, $\pm (1,0,1)$ and $(1,-1,1)$, $(-1,1,1)$, respectively. These solutions correspond to the Frey-curves $E(1,1)$, $E(0,1)$ and $E(1,-1)$ that indeed exist and so there must be newforms associated with them in S1 which a priori will not be possible to eliminate only by comparing the $a_L$.\par 
Going through all the forms in $S1$ and comparing the corresponding values of the ${a_L}'s$ with the possibilities for our Frey-curves we immediately eliminate all except 4 newforms. Here we have eliminated a newform if one of its coefficients $a_L$ is not on the corresponding list above. This can be done because the value of $p$ in the statement of Theorem \ref{grande} is very large hence congruence (\ref{congruencia}), when specified at a trace at $L$ for a prime $L$ of small norm, does not hold modulo such large prime $p$ unless $a_L(f) = a_L(E)$. For example, the first form in the appendix satisfies $a_{L_5}(f) = -9$ and since $a_{L_{5}}(E) \in \{-6, -2 ,2\}$ it is clear that $-9 \equiv -6,-2,2$ (mod $p$) can not hold for $p > 11$. The four remaining newforms correspond to the trivial solutions above plus the twist by $-1$ of $E(1,1)$. The one associated with $E(1,1)$ has level $\mathfrak{P}_{2}^{4}\mathfrak{P}_{13}^{2}$ and the other three $\mathfrak{P}_{2}^{3}\mathfrak{P}_{13}^{2}$. In table \ref{sobras} we list their first eigenvalues.
\begin{table}[ht]
\begin{center}
\begin{tabular}{|c|c|c|c|c|c|c|c|c|c|c|c|}
\hline
		&$a_{L_{3}^{0}}$& $a_{L_{3}^{1}}$& $a_{L_{17}^{0}}$& $a_{L_{17}^{1}}$& $a_{L_{23}^{0}}$& $a_{L_{23}^{1}}$& $a_{L_{5}}$& $a_{L_{29}^{0}}$& $a_{L_{29}^{1}}$& $a_{L_{7}}$& $a_{L_{11}}$\\ \hline

$f_1$&-1&1&7&3&1&7&2&-7&-3&-1&3\\ \hline 
$f_2$&-1&1&3&7&-7&-1&2&-3&-7&-1&3\\ \hline 
$f_3$&-1&-3&-1&-5&5&-9&-6&-3&1&-5&15\\ \hline 
$f_4$&-3&-1&1&-3&-3&-9&-2&-7&5&-11&-15\\ \hline 
 
\end{tabular}
\caption{Values of $a_L$}
\label{sobras}
\end{center}
\end{table} 

To be able to eliminate these newforms we need to use the extra conditions on $d$ and $a+b$. Recall that the solutions $(a,b,c)$ to equation (\ref{caso1}) or (\ref{caso13}) satisfy $d \mid a+b$. Recomputing the possibilities for some $a_L$ but with this extra condition we find that $a_{L_{3}^{0}} = -3$ and $a_{L_{3}^{1}} = -1$ (if $d=3$) , $a_{L_{5}} = -2$ (if $d=5$) , $a_{L_{7}} = -11$ (if $d=7$) or $a_{L_{11}} = -15$ (if $d=11$). By checking in table \ref{sobras} we see that any of the previous conditions is enough to eliminate all $f_i$ except for $f_4$. Actually, $f_4$ is the newform associated with the trivial solution $(1,-1,0)$ and can not be eliminated this way as expected. Finally, if we assume that the solution is first case, Proposition \ref{condE} together with condition $13 \nmid a+b$ guarantees that when restricted to $\mbox{Gal}(\bar{\Q}/K)$ the representations $\rho_{f_4,p}$ and $\rho_{E,p}$ will have different inertia at $\mathfrak{P}_{13}$ and thus can not be isomorphic modulo $\mathfrak{P}$.\\  

To finish the argument we have to eliminate also the newforms in $S2$. Recall that we know the factorization (appendix B) of the characteristic polynomial of $T_{L_3^{0}}$ which we denote by $p_3$. If for $f$ in $S2$ congruence (\ref{congruencia}) holds we also have $$a_{L_3^{0}}(E) \equiv c_{L_3^{0}}(f) \mbox{ (mod } \mathfrak{P}).$$ 
Let $p_c(x)$ be the minimal polynomial of $c_{L_3^{0}}(f)$ which must be a non-linear factor of $p_3$. Thus,
\begin{equation}
\label{congFinal}
p_c(a_{L_3^{0}}(E)) \equiv p_c(c_{L_3^{0}}(f)) \equiv 0 \mbox{ (mod } \mathfrak{P})  
\end{equation}
and $p_c(a_{L_3^{0}}(E)) \neq 0$ because $c_{L_3^{0}}(f) \not\in \mathbb{Z}$. Since $a_{L_3^{0}}(E) \in \{3,-1\}$ by computing $p_c(3)$ and $p_c(-1)$ for all $p_c$ a non-linear factors of $p_3$ we have all the possibilities for $p_c(a_{L_3^{0}}(E))$ and we can see that congruence (\ref{congFinal}) can not hold if $p > 4992539$. Therefore (\ref{congruencia}) also can not hold if $p > 4992539$ and this ends the proof of part (I) in Theorem \ref{grande}. \qed

\begin{rem} It is also possible to eliminate the newforms in $S2$ without knowing the factorization of $p_3$ but this would result in the bound $p > 2^{14546}$ for the exponent. Indeed, let $p_c = \sum r_n x^n$ be the minimal polynomial of a non integer $c_{L_3}(f)$. All the roots $c_{L_3}^\sigma$ satisfy the Weil bound since they are coefficients of the conjugated form $f^\sigma$. Moreover, by knowing the dimension of $S_2(\mathfrak{P}_{2}^s\mathfrak{P}_{13}^{2})$ we can bound all $|r_n|$ using the binomial coefficients. Putting these bounds together we find only a finite number of possibilities for the non-zero value $p_c(a_{L_3}(E))$. The details for this argument can be found in the first version of this work at http://arxiv.org/abs/1112.4521
\end{rem}

Part (II) now follows easily from the proof of part (I). First note that as before it follows from the factorization
\begin{equation}
x^{26} + y^{26} = (x^2 + y^2)\phi(x^2,y^2) = 10 z^p,
\label{eq2} 
\end{equation}
Proposition \ref{trezz2} and Corollary \ref{trezz}, that a solution $(a,b,c)$ must verify $10 \mid a^2 + b^2$.  To a primitive solution $(a,b,c)$ we now attach the Frey-curve $E(a^2,b^2)$. Observe also that $4 \nmid a^2 + b^2$ by looking modulo 4. It now follows from Proposition \ref{condEgamma}, modularity and lowering the level that the set S1 will only have newforms of level $\mathfrak{P}_{2}^4\mathfrak{P}_{13}^{2}$. This means that after comparing the values $a_L$ as before, we eliminate all newforms except for the one corresponding to the curve $E(1,1)$. As we already know, the extra restriction $5 \mid a^2+b^2$ is enough to deal with this newform. In fact, recall that in this case the Frey curve has $a_5 = -2$, and this is different from the corresponding coefficient $a_5$ of $E(1,1)$. The newforms in S2 can be eliminated exactly as in the proof of part (I).\qed  

\newpage

\section{Appendix} 

\subsection{A: Tables with values $a_L$}

\begin{table}[th]
\begin{center}
\begin{tabular}{|c|c|c|c|c|c|c|c|c|c|c|}
\hline
$a_{L_{3}^{0}}$& $a_{L_{3}^{1}}$& $a_{L_{17}^{0}}$& $a_{L_{17}^{1}}$& $a_{L_{23}^{0}}$& $a_{L_{23}^{1}}$& $a_{L_{5}}$& $a_{L_{29}^{0}}$& $a_{L_{29}^{1}}$& $a_{L_{7}}$& $a_{L_{11}}$\\ \hline
-3& -3& -3& -3& -4& -4& -9& 2& 2& -13& -18 \\ \hline
-3&-1& -5& -1& -9& 5& -6& 1& -3& -5& 15\\ \hline 
-3& -1& -1& 3& 3& 9& 2& -7& 5& 11& 15\\ \hline
-3& -1& 1& -3& -3& -9& -2& -7& 5& -11&-15\\ \hline 
-3& -1& 5& 1& 9& -5& 6& 1& -3& 5& -15\\ \hline
-3& 1& -7& -7& 3& -1& 6& 7& -9& 1& 9\\ \hline
-3& 1& 7& 7&-3& 1& -6& 7& -9& -1& -9\\ \hline
-1& -3& -3& 1& -9& -3& -2& 5& -7&-11& -15\\ \hline
-1& -3& -1& -5& 5& -9& -6& -3& 1& -5& 15\\ \hline
-1& -3& 1& 5& -5& 9& 6& -3& 1& 5& -15\\ \hline
-1& -3& 3& -1&9& 3& 2& 5& -7& 11& 15\\ \hline
-1& -1& 3& 3& -6& -6& 1& 0& 0& 5&22\\ \hline
-1& 1& -7& -3& -1& -7& -2& -7& -3& 1& -3\\ \hline
-1& 1& -5& 7& 5& 3& -6& 1& 5& -1& 3\\ \hline
-1& 1& 5& -7& -5& -3& 6& 1& 5& 1& -3\\ \hline
-1& 1& 7& 3& 1& 7& 2& -7& -3& -1& 3\\ \hline
-1& 3& -1& 7& 2& 2& -7& -8& 0& 1& 6\\ \hline
-1& 3& 1& -7& -2& -2& 7& -8& 0& -1& -6\\ \hline
0& 0& 6& 6& 8& 8& -6& 2& 2& -10& -18\\ \hline
1& -3& -7& -7& -1& 3& 6& -9& 7& 1& 9\\ \hline
1& -3& 7& 7& 1& -3& -6& -9& 7& -1& -9\\ \hline
1& -1& -7& 5& -3& -5& 6& 5& 1& 1& -3\\ \hline
1& -1& -3& -7& -7& -1& -2& -3& -7& 1& -3\\ \hline
1& -1& 3& 7& 7& 1& 2& -3& -7& -1& 3\\ \hline
1& -1& 7& -5& 3& 5& -6& 5& 1& -1& 3\\ \hline
1& 1& -3& -3& -1& -1& 6& 3& 3& 13& 21\\ \hline
1& 1& -3& -3& 0& 0& -1& 6& 6& -13& 14\\ \hline
1& 1& -3& -3& 4& 4& -9& -6& -6& 11& -18\\ \hline
1& 1& 3& 3& -4& -4& 9& -6& -6& -11& 18\\ \hline
\end{tabular}
\caption{$a_L$ values for newforms of level $\mathfrak{P}_{2}^{3}\mathfrak{P}_{13}^{2}$}
\label{novas23}
\end{center}
\end{table} 

\newpage

\begin{table}[ht]
\begin{center}
\begin{tabular}{|c|c|c|c|c|c|c|c|c|c|c|}
\hline
$a_{L_{3}^{0}}$& $a_{L_{3}^{1}}$& $a_{L_{17}^{0}}$& $a_{L_{17}^{1}}$& $a_{L_{23}^{0}}$& $a_{L_{23}^{1}}$& $a_{L_{5}}$& $a_{L_{29}^{0}}$& $a_{L_{29}^{1}}$& $a_{L_{7}}$& $a_{L_{11}}$\\ \hline
1& 1& 3& 3& 1& 1& -6& 3& 3& -13& -21\\ \hline
3& -1& -7& 1& -2& -2& 7& 0& -8& -1& -6\\ \hline
3& -1& 7& -1& 2& 2& -7& 0& -8& 1& 6\\ \hline
-3& -3& -3& -3& -1& -1& 6& -1& -1& 13& -3\\ \hline
-3& -3& 3& 3& 1& 1& -6& -1& -1& -13& 3\\ \hline
-3& -1& -7& -7& -3& -1& -6& -7& 9& 1& 9\\ \hline 
-3& -1& -7& 1& -2& 2& -7& 0& 8& -1& -6\\ \hline
-3& -1& 7& -1& 2& -2& 7& 0& 8& 1& 6\\ \hline 
-3& -1& 7& 7& 3& 1& 6& -7& 9& -1& -9\\ \hline 
-3& 1& -7& 1& 2& 2& 7& 0& -8& -1& -6\\ \hline 
-3& 1& -5& -1& 9& 5& 6& -1& 3& -5& 15\\ \hline
 -3& 1& -1& 3& -3& 9& -2& 7& -5& 11& 15\\ \hline 
-3& 1& 1& -3& 3& -9& 2& 7& -5& -11& -15\\ \hline 
-3& 1& 5& 1& -9& -5& -6& -1& 3& 5& -15\\ \hline
 -3& 1& 7& -1& -2& -2& -7& 0& -8& 1& 6\\ \hline
-3& 3& -3& -3& 1& -1& -6& 1& 1& 13& -3\\ \hline
 -3& 3& -3& -3& 4& -4& 9& -2& -2& -13& -18\\ \hline
 -3& 3& 3& 3& -4& 4& -9& -2& -2& 13& 18\\ \hline 
-3& 3& 3& 3&-1& 1& 6& 1& 1& -13& 3\\ \hline 
-2& -2& -3& -3& 6& 6& 7& 3& 3& 14& 22\\ \hline
 -2&-2& 3& 3& -6& -6& -7& 3& 3& -14& -22\\ \hline 
-2& 2& -3& -3& -6& 6& -7& -3&-3& 14& 22\\ \hline 
-2& 2& -3& -3& -6& 6& -1& -9& -9& -2& 22\\ \hline 
-2& 2& 3& 3&6& -6& 1& -9& -9& 2& -22\\ \hline 
-2& 2& 3& 3& 6& -6& 7& -3& -3& -14& -22\\ \hline
-1& -3& -7& -7& -1& -3& -6& 9& -7& 1& 9\\ \hline 
-1& -3& -1& 7& -2& 2& 7& 8&0& 1& 6\\ \hline
 -1& -3& 1& -7& 2& -2& -7& 8& 0& -1& -6\\ \hline
 -1& -3& 7& 7& 1& 3& 6& 9& -7& -1& -9\\ \hline
 -1& -1& -7& -3& 1& -7& 2& 7& 3& 1& -3\\ \hline 
-1& -1& -7& 5& -3& 5& -6& -5& -1& 1& -3\\ \hline
 -1& -1& -5& 7& -5& 3& 6& -1& -5& -1& 3\\ \hline
 -1& -1& -3& -7& -7& 1& 2& 3& 7& 1& -3\\ \hline
\end{tabular}
\caption{$a_L$ values for newforms of level $\mathfrak{P}_{2}^{4}\mathfrak{P}_{13}^{2}$}
\label{novas241}
\end{center}
\end{table} 

\newpage

\begin{table}[ht]
\begin{center}
\begin{tabular}{|c|c|c|c|c|c|c|c|c|c|c|}
\hline
$a_{L_{3}^{0}}$& $a_{L_{3}^{1}}$& $a_{L_{17}^{0}}$& $a_{L_{17}^{1}}$& $a_{L_{23}^{0}}$& $a_{L_{23}^{1}}$& $a_{L_{5}}$& $a_{L_{29}^{0}}$& $a_{L_{29}^{1}}$& $a_{L_{7}}$& $a_{L_{11}}$\\ \hline
-1& -1& -3& -3& 1& 1& 6& 3& 3& 13& 21\\ \hline
-1& -1& 3& 3& -1& -1& -6& 3& 3& -13& -21\\ \hline
-1& -1& 3& 3& 0& 0& 1& 6& 6& 13& -14\\ \hline 
-1& -1& 3& 3& 4& 4& 9& -6& -6& -11& 18\\ \hline 
-1& -1& 3& 7& 7& -1& -2& 3& 7& -1& 3\\ \hline 
-1& -1& 5& -7& 5& -3& -6& -1& -5& 1& -3\\ \hline 
-1& -1& 7& -5& 3& -5& 6& -5& -1& -1& 3\\ \hline
-1& -1& 7& 3& -1& 7& -2& 7& 3& -1& 3\\ \hline 
-1& 1& -7& 5& 3& 5& 6& 5& 1& 1& -3\\ \hline 
-1& 1& -3& -7& 7& 1& -2& -3& -7& 1& -3\\ \hline
-1& 1& -3& -3& -6& 6& 1& 0& 0& -5& -22\\ \hline
-1& 1& -3& -3& -1& 1& -6& -3& -3& 13& 21\\ \hline
-1& 1& -3&-3& 0& 0& 1& -6& -6& -13& 14\\ \hline
-1& 1& -3& -3& 3& -3& 10& 9& 9& -11& 5\\ \hline
-1& 1& -3& -3& 4& -4& 9& 6& 6& 11& -18\\ \hline 
-1& 1& 3& 3& -4& 4& -9& 6& 6& -11& 18\\ \hline
-1& 1& 3& 3& -3& 3& -10& 9& 9& 11& -5\\ \hline
-1& 1& 3& 3& 0& 0& -1& -6& -6& 13& -14\\ \hline 
-1& 1& 3& 3& 1& -1& 6& -3& -3& -13& -21\\ \hline
-1& 1&3& 3& 6& -6& -1& 0& 0& 5& 22\\ \hline 
-1& 1& 3& 7& -7& -1& 2& -3& -7& -1& 3\\ \hline
-1& 1& 7& -5& -3& -5& -6& 5& 1& -1& 3\\ \hline 
-1& 3& -7& -7& 1& -3& 6& -9& 7& 1& 9\\ \hline 
-1& 3& -3& 1& 9& -3& 2& -5& 7& -11& -15\\ \hline 
-1& 3& -1& -5& -5& -9& 6& 3& -1& -5& 15\\ \hline
-1& 3& 1& 5& 5& 9& -6& 3& -1& 5& -15\\ \hline 
-1& 3& 3& -1& -9& 3& -2& -5& 7& 11& 15\\ \hline 
-1& 3& 7& 7& -1& 3& -6& -9& 7& -1& -9\\ \hline 
0& 0& -6& -6& -8& 8& -6& -2& -2& 10& 18\\ \hline 
0& 0& -6& -6& 8& -8& -6& -2& -2& 10& 18\\ \hline 
0& 0& -6& -6& 8& 8& 6& 2& 2& 10& 18\\ \hline 
0& 0& -3& -3& -4& -4& 9& -1& -1& -2& 6\\ \hline 
\end{tabular}
\caption{$a_L$ values for newforms of level $\mathfrak{P}_{2}^{4}\mathfrak{P}_{13}^{2}$ (cont.)}
\label{novas242}
\end{center}
\end{table} 

\newpage

\begin{table}[ht]
\begin{center}
\begin{tabular}{|c|c|c|c|c|c|c|c|c|c|c|}
\hline
$a_{L_{3}^{0}}$& $a_{L_{3}^{1}}$& $a_{L_{17}^{0}}$& $a_{L_{17}^{1}}$& $a_{L_{23}^{0}}$& $a_{L_{23}^{1}}$& $a_{L_{5}}$& $a_{L_{29}^{0}}$& $a_{L_{29}^{1}}$& $a_{L_{7}}$& $a_{L_{11}}$\\ \hline
0& 0& -3& -3& -4& 4& -9& 1& 1& -2& 6\\ \hline
0& 0& -3& -3& 4& -4& -9& 1& 1& -2& 6\\ \hline 
0& 0& 3& 3& -4& 4& 9& 1& 1& 2& -6\\ \hline 
0& 0& 3& 3& 4& -4& 9& 1& 1& 2& -6\\ \hline 
0& 0& 3& 3& 4& 4& -9& -1& -1& 2& -6\\ \hline 
0& 0& 6& 6& -8& -8& -6& 2& 2& -10& -18\\ \hline 
0& 0& 6& 6& -8& 8& 6& -2& -2& -10& -18\\ \hline 
0& 0& 6& 6& 8& -8& 6& -2& -2& -10& -18\\ \hline 
1& -3& -3& 1& -9& 3& 2& -5& 7& -11& -15\\ \hline 
1& -3& -1& -5& 5& 9& 6& 3& -1& -5& 15\\ \hline
1& -3& -1& 7& -2& -2& -7& -8& 0& 1& 6\\ \hline 
1& -3& 1& -7& 2& 2& 7& -8& 0& -1& -6\\ \hline 
1& -3& 1& 5& -5& -9& -6& 3& -1& 5& -15\\ \hline 
1& -3& 3& -1& 9& -3& -2& -5& 7& 11& 15\\ \hline 
1& -1& -7& -3& 1& 7& -2& -7& -3& 1& -3\\ \hline 
1& -1& -5& 7& -5& -3& -6& 1& 5& -1& 3\\ \hline 
1& -1& -3& -3& -4& 4& 9& 6& 6& 11& -18\\ \hline 
1& -1& -3& -3& -3& 3& 10& 9& 9& -11& 5\\ \hline 
1& -1& -3& -3& 0& 0& 1& -6& -6& -13& 14\\ \hline 
1& -1& -3& -3& 1& -1& -6& -3& -3& 13& 21\\ \hline 
1& -1& -3& -3& 6& -6& 1& 0& 0& -5& -22\\ \hline 
1& -1& 3& 3& -6& 6& -1& 0& 0& 5& 22\\ \hline
1& -1& 3& 3& -1& 1& 6& -3& -3& -13& -21\\ \hline 
1& -1& 3& 3& 0& 0& -1& -6& -6& 13& -14\\ \hline 
1& -1& 3& 3& 3& -3& -10& 9& 9& 11& -5\\ \hline 
1& -1& 3& 3& 4& -4& -9& 6& 6& -11& 18\\ \hline 
1& -1& 5& -7& 5& 3& 6& 1& 5& 1& -3\\ \hline 
1& -1& 7& 3& -1& -7& 2& -7& -3& -1& 3\\ \hline 
1& 1& -7& -3& -1& 7& 2& 7& 3& 1& -3\\ \hline 
1& 1& -7& 5& 3& -5& -6& -5& -1& 1& -3\\ \hline 
1& 1& -5& 7& 5& -3& 6& -1& -5& -1& 3\\ \hline 
1& 1& -3& -7& 7& -1& 2& 3& 7& 1& -3\\ \hline 
\end{tabular}
\caption{$a_L$ values for newforms of level $\mathfrak{P}_{2}^{4}\mathfrak{P}_{13}^{2}$ (cont.)}
\label{novas243}
\end{center}
\end{table} 

\newpage

\begin{table}[ht]
\begin{center}
\begin{tabular}{|c|c|c|c|c|c|c|c|c|c|c|}
\hline
$a_{L_{3}^{0}}$& $a_{L_{3}^{1}}$& $a_{L_{17}^{0}}$& $a_{L_{17}^{1}}$& $a_{L_{23}^{0}}$& $a_{L_{23}^{1}}$& $a_{L_{5}}$& $a_{L_{29}^{0}}$& $a_{L_{29}^{1}}$& $a_{L_{7}}$& $a_{L_{11}}$\\ \hline
1& 1& -3& -3& -6& -6& -1& 0& 0& -5& -22\\ \hline 
1& 1& -3& -3& 3& 3& -10& -9& -9& -11& 5\\ \hline 
1& 1& 3& 3& -3& -3& 10& -9& -9& 11& -5\\ \hline 
1& 1& 3& 3& 6& 6& 1& 0& 0& 5& 22\\ \hline 
1& 1& 3& 7& -7& 1& -2& 3& 7& -1& 3\\ \hline 
1& 1& 5& -7& -5& 3& -6& -1& -5& 1& -3\\ \hline 
1& 1& 7& -5& -3& 5& 6& -5& -1& -1& 3\\ \hline 
1& 1& 7& 3& 1& -7& -2& 7& 3& -1& 3\\ \hline 
1& 3& -7& -7& 1& 3& -6& 9& -7& 1& 9\\ \hline 
1& 3& -3& 1& 9& 3& -2& 5& -7& -11& -15\\ \hline 
1& 3& -1& -5& -5& 9& -6& -3& 1& -5& 15\\ \hline 
1& 3& -1& 7& 2& -2& 7& 8& 0& 1& 6\\ \hline 
1& 3& 1& -7& -2& 2& -7& 8& 0& -1& -6\\ \hline 
1& 3& 1& 5& 5& -9& 6& -3& 1& 5& -15\\ \hline 
1& 3& 3& -1& -9& -3& 2& 5& -7& 11& 15\\ \hline
1& 3& 7& 7& -1& -3& 6& 9& -7& -1& -9\\ \hline 
2& -2& -3& -3& 6& -6& -7& -3& -3& 14& 22\\ \hline 
2& -2& -3& -3& 6& -6& -1& -9& -9& -2& 22\\ \hline 
2& -2& 3& 3& -6& 6& 1& -9& -9& 2& -22\\ \hline 
2& -2& 3& 3& -6& 6& 7& -3& -3& -14& -22\\ \hline
2& 2& -3& -3& -6& -6& 1& 9& 9& -2& 22\\ \hline 
2& 2& 3& 3& 6& 6& -1& 9& 9& 2& -22\\ \hline 
3& -3& -3& -3& -4& 4& 9& -2& -2& -13& -18\\ \hline 
3& -3& -3& -3& -1& 1& -6& 1& 1& 13& -3\\ \hline 
3& -3& 3& 3& 1& -1& 6& 1& 1& -13& 3\\ \hline 
3& -3& 3& 3& 4& -4& -9& -2& -2& 13& 18\\ \hline 
3& -1& -7& -7& -3& 1& 6& 7& -9& 1& 9\\ \hline
3& -1& -5& -1& -9& -5& 6& -1& 3& -5& 15\\ \hline 
3& -1& -1& 3& 3& -9& -2& 7& -5& 11& 15\\ \hline 
3& -1& 1& -3& -3& 9& 2& 7& -5& -11& -15\\ \hline 
3& -1& 5& 1& 9& 5& -6& -1& 3& 5& -15\\ \hline 
\end{tabular}
\caption{$a_L$ values for newforms of level $\mathfrak{P}_{2}^{4}\mathfrak{P}_{13}^{2}$ (cont.)}
\label{novas244}
\end{center}
\end{table} 

\newpage

\begin{table}[ht]
\begin{center}
\begin{tabular}{|c|c|c|c|c|c|c|c|c|c|c|}
\hline
$a_{L_{3}^{0}}$& $a_{L_{3}^{1}}$& $a_{L_{17}^{0}}$& $a_{L_{17}^{1}}$& $a_{L_{23}^{0}}$& $a_{L_{23}^{1}}$& $a_{L_{5}}$& $a_{L_{29}^{0}}$& $a_{L_{29}^{1}}$& $a_{L_{7}}$& $a_{L_{11}}$\\ \hline
3& -1& 7& 7& 3& -1& -6& 7& -9& -1& -9\\ \hline 
3& 1& -7& -7& 3& 1& -6& -7& 9& 1& 9\\ \hline 
3& 1& -7& 1& 2& -2& -7& 0& 8& -1& -6\\ \hline
3& 1& -5& -1& 9& -5& -6& 1& -3& -5& 15\\ \hline 
3& 1& -1& 3& -3& -9& 2& -7& 5& 11& 15\\ \hline 
3& 1& 1& -3& 3& 9& -2& -7& 5& -11& -15\\ \hline 
3& 1& 5& 1& -9& 5& 6& 1& -3& 5& -15\\ \hline 
3& 1& 7& -1& -2& 2& 7& 0& 8& 1& 6\\ \hline
3& 1& 7& 7& -3& -1& 6& -7& 9& -1& -9\\ \hline 
3& 3& -3& -3& 4& 4& -9& 2& 2& -13& -18\\ \hline 
3& 3& 3& 3& -4& -4& 9& 2& 2& 13& 18 \\ \hline
-1& -1& -3& -3& -4& -4& -9& -6& -6& 11& -18\\ \hline
-1& -1& -3& -3& 0& 0& -1& 6& 6& -13& 14\\ \hline 
\end{tabular}
\caption{$a_L$ values for newforms of level $\mathfrak{P}_{2}^{4}\mathfrak{P}_{13}^{2}$ (cont.)}
\label{novas245}
\end{center}
\end{table}

\subsection{B: Factorization of $p_3$}

The polynomial $p_3$ on $S_2(\mathfrak{P}_{2}^{4}\mathfrak{P}_{13}^{2})$ has the following factors:\\

$[x - 3, 38],[x - 2, 16],[x - 1, 84],[x, 32],[x + 1, 92],[x + 2, 12],[x + 3, 45][x^2 - 4 x + 2, 6],[x^2 - 3 x - 1, 17],[x^2 - 3 x + 1, 39],[x^2 - 2 x - 2, 4],[x^2 - 2 x - 1, 10],[x^2 - x - 5, 18],[x^2 - x - 4, 23],[x^2 - x - 3, 33],[x^2 - x - 1, 25],[x^2 - 8, 10],[x^2 - 5, 12],[x^2 - 3, 20],[x^2 - 2, 10],[x^2 + x - 5, 6],[x^2 + x - 4, 22],[x^2 + x - 3, 42],[x^2 + x - 1, 30],[x^2 + 2 x - 2, 10],[x^2 + 2 x - 1, 4],[x^2 + 3 x - 1, 21],[x^2 + 3 x + 1, 37],[x^2 + 4 x + 2, 4],[x^3 - 4 x^2 + 3 x + 1, 4],[x^3 - 3 x^2 - 4 x + 13, 10],[x^3 - 2 x^2 - 4 x + 7, 6],[x^3 - 2 x^2 - x + 1, 4],[x^3 - 7 x - 7, 4],[x^3 - 7 x + 7, 8],[x^3 - 4 x - 1, 9],[x^3 - 4 x + 1, 6],[x^3 + 2 x^2 - 4 x - 7, 9],[x^3 + 2 x^2 - x - 1, 12],[x^3 + 3 x^2 - 4 x - 13, 4],[x^3 + 4 x^2 + 3 x - 1, 6],[x^4 - 4 x^3 - x^2 + 10 x + 2, 6],[x^4 - 4 x^3 + 8 x - 1, 4],  [x^4 - 3 x^3 - 6 x^2 + 23 x - 13, 4],[x^4 - 2 x^3 - 9 x^2 + 22 x - 11, 12],[x^4 - 2 x^3 - 7 x^2 + 8 x - 1, 6],[x^4 - 2 x^3 - 7 x^2 + 8 x + 4, 4],[x^4 - 2 x^3 - 6 x^2 + 10 x + 1, 14],[x^4 - 2 x^3 - 4 x^2 + 8 x - 2, 4],[x^4 + 2 x^3 - 9 x^2 - 22 x - 11, 6],[x^4 + 2 x^3 - 7 x^2 - 8 x - 1, 9],[x^4 + 2 x^3 - 7 x^2 - 8 x + 4, 6],[x^4 + 2 x^3 - 6 x^2 - 10 x + 1, 8],[x^4 + 2 x^3 - 4 x^2 - 8 x - 2, 10],[x^4 + 3 x^3 - 6 x^2 - 23 x - 13, 6],[x^4 + 4 x^3 - x^2 - 10 x + 2, 4],[x^4 + 4 x^3 - 8 x - 1, 8],[x^6 - 3 x^5 - 12 x^4 + 36 x^3 + 18 x^2 - 68 x + 29, 10],[x^6 - 3 x^5 - 11 x^4 + 31 x^3 + 15 x^2 - 25 x - 4, 9],[x^6 - 3 x^5 - 9 x^4 + 23 x^3 + 15 x^2 - 13 x - 1, 6],[x^6 - 2 x^5 - 11 x^4 + 16 x^3 + 35 x^2 - 26 x - 25, 4],[x^6 - x^5 - 13 x^4 + 11 x^3 + 49 x^2 - 27 x - 52, 9],[x^6 + x^5 - 13 x^4 - 11 x^3 + 49 x^2 + 27 x - 52, 6],[x^6 + 2 x^5 - 11 x^4 - 16 x^3 + 35 x^2 + 26 x - 25, 8],[x^6 + 3 x^5 - 12 x^4 - 36 x^3 + 18 x^2 + 68 x + 29, 4],[x^6 + 3 x^5 - 11 x^4 - 31 x^3 + 15 x^2 + 25 x - 4, 6],[x^6 + 3 x^5 - 9 x^4 - 23 x^3 + 15 x^2 + 13 x - 1, 4],[x^8 - 6 x^7 + x^6 + 48 x^5 - 65 x^4 - 54 x^3 + 115 x^2 - 48 x + 4, 8],[x^8 - 3 x^7 - 20 x^6 + 57 x^5 + 124 x^4 - 327 x^3 - 245 x^2 + 588 x + 16,4],[x^8 + 3 x^7 - 20 x^6 - 57 x^5 + 124 x^4 + 327 x^3 - 245 x^2 - 588 x + 16,2],[x^8 + 6 x^7 + x^6 - 48 x^5 - 65 x^4 + 54 x^3 + 115 x^2 + 48 x + 4, 4],[x^9 - 4 x^8 - 9 x^7 + 50 x^6 - 5 x^5 - 156 x^4 + 125 x^3 + 50 x^2 - 40 x +4, 6],[x^9 - 3 x^8 - 11 x^7 + 32 x^6 + 38 x^5 - 100 x^4 - 47 x^3 + 75 x^2 + 37 x +4, 4],[x^9 - 2 x^8 - 17 x^7 + 34 x^6 + 75 x^5 - 158 x^4 - 31 x^3 + 106 x^2 - 20 x- 4, 4],[x^9 - x^8 - 19 x^7 + 16 x^6 + 114 x^5 - 76 x^4 - 251 x^3 + 165 x^2 + 181 x- 128, 6],[x^9 + x^8 - 19 x^7 - 16 x^6 + 114 x^5 + 76 x^4 - 251 x^3 - 165 x^2 + 181 x+ 128, 4],[x^9 + 2 x^8 - 17 x^7 - 34 x^6 + 75 x^5 + 158 x^4 - 31 x^3 - 106 x^2 - 20 x+ 4, 6],[x^9 + 3 x^8 - 11 x^7 - 32 x^6 + 38 x^5 + 100 x^4 - 47 x^3 - 75 x^2 + 37 x -4, 6],[x^9 + 4 x^8 - 9 x^7 - 50 x^6 - 5 x^5 + 156 x^4 + 125 x^3 - 50 x^2 - 40 x -4, 4],[x^{12} - x^{11} - 22 x^{10} + 30 x^9 + 153 x^8 - 276 x^7 - 317 x^6 + 863 x^5 -182 x^4 - 513 x^3 + 242 x^2 + 22 x - 1, 4],[x^{12} + x^{11} - 22 x^{10} - 30 x^9 + 153 x^8 + 276 x^7 - 317 x^6 - 863 x^5 -182 x^4 + 513 x^3 + 242 x^2 - 22 x - 1, 8],[x^{16} - 2 x^{15} - 26 x^{14} + 48 x^{13} + 261 x^{12} - 442 x^{11} - 1300 x^{10} +2024 x^9 + 3449 x^8 - 4958 x^7 - 4874 x^6 + 6520 x^5 + 3355 x^4 -4294 x^3 - 776 x^2 + 1104 x - 74, 6],[x^{16} + 2 x^{15} - 26 x^{14} - 48 x^{13} + 261 x^{12} + 442 x^{11} - 1300 x^{10} -2024 x^9 + 3449 x^8 + 4958 x^7 - 4874 x^6 - 6520 x^5 + 3355 x^4 +4294 x^3 - 776 x^2 - 1104 x - 74, 4],[x^{18} - 6 x^{17} - 22 x^{16} + 184 x^{15} + 101 x^{14} - 2206 x^{13} + 1048 x^{12} +13080 x^{11} - 12983 x^{10} - 39490 x^9 + 52906 x^8 + 54352 x^7 - 91701 x^6- 23122 x^5 + 56412 x^4 + 5600 x^3 - 13034 x^2 - 1480 x + 568, 4],[x^{18} - 4 x^{17} - 33 x^{16} + 150 x^{15} + 387 x^{14} - 2233 x^{13} - 1578 x^{12} +16799 x^{11} - 4197 x^{10} - 65971 x^9 + 59322 x^8 + 117310 x^7 - 188308 x^6- 21264 x^5 + 190984 x^4 - 132221 x^3 + 32604 x^2 - 1200 x - 379, 6],[x^{18} - 3 x^{17} - 28 x^{16} + 82 x^{15} + 318 x^{14} - 915 x^{13} - 1870 x^{12} +5426 x^{11} + 5939 x^{10} - 18603 x^9 - 9152 x^8 + 37297 x^7 + 2528 x^6 -41228 x^5 + 10028 x^4 + 20669 x^3 - 9779 x^2 - 1970 x + 1259, 4],[x^{18} - 3 x^{17} - 28 x^{16} + 84 x^{15} + 308 x^{14} - 921 x^{13} - 1692 x^{12} +4994 x^{11} + 4927 x^{10} - 13943 x^9 - 7648 x^8 + 19011 x^7 + 6382 x^6 -10700 x^5 - 3212 x^4 + 2003 x^3 + 627 x^2 - 20 x - 1, 4],[x^{18} - 3 x^{17} - 26 x^{16} + 78 x^{15} + 270 x^{14} - 815 x^{13} - 1436 x^{12} +4416 x^{11} + 4153 x^{10} - 13377 x^9 - 6258 x^8 + 22693 x^7 + 3772 x^6 -20184 x^5 + 782 x^4 + 7699 x^3 - 1277 x^2 - 514 x + 13, 4],[x^{18} - 3 x^{17} - 26 x^{16} + 80 x^{15} + 264 x^{14} - 837 x^{13} - 1362 x^{12} +4500 x^{11} + 3737 x^{10} - 13433 x^9 - 4682 x^8 + 22123 x^7 - 202 x^6 -18304 x^5 + 5878 x^4 + 5249 x^3 - 3531 x^2 + 688 x - 43, 4],[x^{18} - 2 x^{17} - 34 x^{16} + 64 x^{15} + 461 x^{14} - 802 x^{13} - 3208 x^{12} +      4968 x^{11} + 12385 x^{10} - 15862 x^9 - 26994 x^8 + 24760 x^7 + 31587 x^6 -14558 x^5 - 16140 x^4 - 640 x^3 + 910 x^2 + 56 x - 8, 6],[x^{18} + 2 x^{17} - 34 x^{16} - 64 x^{15} + 461 x^{14} + 802 x^{13} - 3208 x^{12} -4968 x^{11} + 12385 x^{10} + 15862 x^9 - 26994 x^8 - 24760 x^7 + 31587 x^6 +14558 x^5 - 16140 x^4 + 640 x^3 + 910 x^2 - 56 x - 8, 4],[x^{18} + 3 x^{17} - 28 x^{16} - 84 x^{15} + 308 x^{14} + 921 x^{13} - 1692 x^{12} -4994 x^{11} + 4927 x^{10} + 13943 x^9 - 7648 x^8 - 19011 x^7 + 6382 x^6 +10700 x^5 - 3212 x^4 - 2003 x^3 + 627 x^2 + 20 x - 1, 6],[x^{18} + 3 x^{17} - 28 x^{16} - 82 x^{15} + 318 x^{14} + 915 x^{13} - 1870 x^{12} -5426 x^{11} + 5939 x^{10} + 18603 x^9 - 9152 x^8 - 37297 x^7 + 2528 x^6 +41228 x^5 + 10028 x^4 - 20669 x^3 - 9779 x^2 + 1970 x + 1259, 8],[x^{18} + 3 x^{17} - 26 x^{16} - 80 x^{15} + 264 x^{14} + 837 x^{13} - 1362 x^{12} -4500 x^{11} + 3737 x^{10} + 13433 x^9 - 4682 x^8 - 22123 x^7 - 202 x^6 +18304 x^5 + 5878 x^4 - 5249 x^3 - 3531 x^2 - 688 x - 43, 6],[x^{18} + 3 x^{17} - 26 x^{16} - 78 x^{15} + 270 x^{14} + 815 x^{13} - 1436 x^{12} -       4416 x^{11} + 4153 x^{10} + 13377 x^9 - 6258 x^8 - 22693 x^7 + 3772 x^6 +20184 x^5 + 782 x^4 - 7699 x^3 - 1277 x^2 + 514 x + 13, 6], [x^{18} + 4 x^{17} - 33 x^{16} - 150 x^{15} + 387 x^{14} + 2233 x^{13} - 1578 x^{12} -16799 x^{11} - 4197 x^{10} + 65971 x^9 + 59322 x^8 - 117310 x^7 - 188308 x^6+ 21264 x^5 + 190984 x^4 + 132221 x^3 + 32604 x^2 + 1200 x - 379, 4],[x^{18} + 6 x^{17} - 22 x^{16} - 184 x^{15} + 101 x^{14} + 2206 x^{13} + 1048 x^{12} -13080 x^{11} - 12983 x^{10} + 39490 x^9 + 52906 x^8 - 54352 x^7 - 91701 x^6+ 23122 x^5 + 56412 x^4 - 5600 x^3 - 13034 x^2 + 1480 x + 568, 6],[x^{27} - 3 x^{26} - 55 x^{25} + 164 x^{24} + 1304 x^{23} - 3858 x^{22} - 17519 x^{21} + 51323 x^{20} + 147499 x^{19} - 427049 x^{18} - 812192 x^{17} + 2322310 x^{16} +2957083 x^{15} - 8374150 x^{14} - 7010721 x^{13} + 19892990 x^{12} + 10323950 x^{11} - 30317349 x^{10} - 8509315 x^9 + 28202264 x^8 + 2963409 x^7- 14792310 x^6 + 220229 x^5 + 3889959 x^4 - 397376 x^3 - 380960 x^2 +72460 x - 2647, 6],[x^{27} - 3 x^{26} - 53 x^{25} + 160 x^{24} + 1210 x^{23} - 3680 x^{22} - 15631 x^{21} +    47937 x^{20} + 126405 x^{19} - 390929 x^{18} - 670268 x^{17} + 2085994 x^{16} +2381831 x^{15} - 7410738 x^{14} - 5717857 x^{13} + 17541450 x^{12} +9224648 x^{11} - 27292133 x^{10} - 9719725 x^9 + 27067090 x^8 + 6187223 x^7- 16180684 x^6 - 1937401 x^5 + 5264635 x^4 + 99800 x^3 - 738836 x^2 +47264 x + 14987, 6],[x^{27} + 3 x^{26} - 55 x^{25} - 164 x^{24} + 1304 x^{23} + 3858 x^{22} - 17519 x^{21} -51323 x^{20} + 147499 x^{19} + 427049 x^{18} - 812192 x^{17} - 2322310 x^{16} +2957083 x^{15} + 8374150 x^{14} - 7010721 x^{13} - 19892990 x^{12} +       10323950 x^{11} + 30317349 x^{10} - 8509315 x^9 - 28202264 x^8 + 2963409 x^7+ 14792310 x^6 + 220229 x^5 - 3889959 x^4 - 397376 x^3 + 380960 x^2 +72460 x + 2647, 4],[x^{27} + 3 x^{26} - 53 x^{25} - 160 x^{24} + 1210 x^{23} + 3680 x^{22} - 15631 x^{21} -47937 x^{20} + 126405 x^{19} + 390929 x^{18} - 670268 x^{17} - 2085994 x^{16} +2381831 x^{15} + 7410738 x^{14} - 5717857 x^{13} - 17541450 x^{12} +       9224648 x^{11} + 27292133 x^{10} - 9719725 x^9 - 27067090 x^8 + 6187223 x^7+ 16180684 x^6 - 1937401 x^5 - 5264635 x^4 + 99800 x^3 + 738836 x^2 +47264 x - 14987, 4]$\\

The polynomial $p_3$ on $S_2(\mathfrak{P}_{2}^{3}\mathfrak{P}_{13}^{2})$ has the following factors:\\

$[x - 3, 8],[x - 2, 8],[x - 1, 20],[x, 12],[x + 1, 28],[x + 2, 4],[x + 3, 15],[x^2 - 4x + 2, 2],[x^2 - 3x - 1, 3],[x^2 - 3x + 1, 11],[x^2 - 2x - 1, 6],
[x^2 - x - 5, 12],[x^2 - x - 4, 5],[x^2 - x - 3, 3],[x^2 - x - 1, 5],[x^2 - 8, 2],[x^2 - 5, 4],[x^2 - 3, 4],[x^2 - 2, 2],[x^2 + x - 4, 4],[x^2 + x - 3, 12],[x^2 + x - 1, 10],[x^2 + 2x - 2, 6],[x^2 + 3x - 1, 7],[x^2 + 3x + 1, 9],[x^3 - 3x^2 - 4x + 13, 6],[x^3 - 7x + 7, 4],[x^3 - 4x - 1, 3],[x^3 + 2x^2 - 4x - 7, 3],[x^3 + 2x^2 - x - 1, 8],[x^3 + 4x^2 + 3x - 1, 2],[x^4 - 4x^3 - x^2 + 10x + 2, 2],[x^4 - 2x^3 - 9x^2 + 22x - 11, 6],[x^4 - 2x^3 - 6x^2 + 10x + 1, 6],[x^4 + 2x^3 - 7x^2 - 8x - 1, 3],[x^4 + 2x^3 - 7x^2 - 8x + 4, 2],[x^4 + 2x^3 - 4x^2 - 8x - 2, 6],[x^4 + 3x^3 - 6x^2 - 23x - 13, 2],
[x^4 + 4x^3 - 8x - 1, 4],[x^6 - 3x^5 - 12x^4 + 36x^3 + 18x^2 - 68x + 29, 6],[x^6 - 3x^5 - 11x^4 + 31x^3 + 15x^2 - 25x - 4, 3],[x^6 - 3x^5 - 9x^4 + 23x^3 + 15x^2 - 13x - 1, 2],[x^6 - x^5 - 13x^4 + 11x^3 + 49x^2 - 27x - 52, 3],[x^6 + 2x^5 - 11x^4 - 16x^3 + 35x^2 + 26x - 25, 4],[x^8 - 6x^7 + x^6 + 48x^5 - 65x^4 - 54x^3 + 115x^2 - 48x + 4, 4],[x^8 - 3x^7 - 20x^6 + 57x^5 + 124x^4 - 327x^3 - 245x^2 + 588x + 16, 2],[x^9 - 4x^8 - 9x^7 + 50x^6 - 5x^5 - 156x^4 + 125x^3 + 50x^2 - 40x +4, 2],[x^9 - x^8 - 19x^7 + 16x^6 + 114x^5 - 76x^4 - 251x^3 + 165x^2 + 181x- 128, 2],[x^9 + 2x^8 - 17x^7 - 34x^6 + 75x^5 + 158x^4 - 31x^3 - 106x^2 - 20x+ 4, 2],[x^9 + 3x^8 - 11x^7 - 32x^6 + 38x^5 + 100x^4 - 47x^3 - 75x^2 + 37x - 4, 2],[x^{12} + x^{11} - 22x^{10} - 30x^9 + 153x^8 + 276x^7 - 317x^6 - 863x^5 -182x^4 + 513x^3 + 242x^2 - 22x - 1, 4],[x^{16} - 2x^{15} - 26x^{14} + 48x^{13} + 261x^{12} - 442x^{11} - 1300x^{10} +   2024x^9 + 3449x^8 - 4958x^7 - 4874x^6 + 6520x^5 + 3355x^4 - 4294x^3 - 776x^2 + 1104x - 74, 2],[x^{18} - 4x^{17} - 33x^{16} + 150x^{15} + 387x^{14} - 2233x^{13} - 1578x^{12} + 16799x^{11} - 4197x^{10} - 65971x^9 + 59322x^8 + 117310x^7 - 188308x^6- 21264x^5 + 190984x^4 - 132221x^3 + 32604x^2 - 1200x - 379, 2],[x^{18} - 2x^{17} - 34x^{16} + 64x^{15} + 461x^{14} - 802x^{13} - 3208x^{12} + 4968x^{11} + 12385x^{10} - 15862x^9 - 26994x^8 + 24760x^7 + 31587x^6 -14558x^5 - 16140x^4 - 640x^3 + 910x^2 + 56x - 8, 2],[x^{18} + 3x^{17} - 28x^{16} - 84x^{15} + 308x^{14} + 921x^{13} - 1692x^{12} -4994x^{11} + 4927x^{10} + 13943x^9 - 7648x^8 - 19011x^7 + 6382x^6 +10700x^5 - 3212x^4 - 2003x^3 + 627x^2 + 20x - 1, 2],[x^{18} + 3x^{17} - 28x^{16} - 82x^{15} + 318x^{14} + 915x^{13} - 1870x^{12} -5426x^{11} + 5939x^{10} + 18603x^9 - 9152x^8 - 37297x^7 + 2528x^6 +41228x^5 + 10028x^4 - 20669x^3 - 9779x^2 + 1970x + 1259, 4],[x^{18} + 3x^{17} - 26x^{16} - 80x^{15} + 264x^{14} + 837x^{13} - 1362x^{12} -4500x^{11} + 3737x^{10} + 13433x^9 - 4682x^8 - 22123x^7 - 202x^6 + 18304x^5 + 5878x^4 - 5249x^3 - 3531x^2 - 688x - 43, 2],[x^{18} + 3x^{17} - 26x^{16} - 78x^{15} + 270x^{14} + 815x^{13} - 1436x^{12} -4416x^{11} + 4153x^{10} + 13377x^9 - 6258x^8 - 22693x^7 + 3772x^6 + 20184x^5 + 782x^4 - 7699x^3 - 1277x^2 + 514x + 13, 2],[x^{18} + 6x^{17} - 22x^{16} - 184x^{15} + 101x^{14} + 2206x^{13} + 1048x^{12} -13080x^{11} - 12983x^{10} + 39490x^9 + 52906x^8 - 54352x^7 - 91701x^6 + 23122x^5 + 56412x^4 - 5600x^3 - 13034x^2 + 1480x + 568, 2],[x^{27} - 3x^{26} - 55x^{25} + 164x^{24} + 1304x^{23} - 3858x^{22} - 17519x^{21} + 51323x^{20} + 147499x^{19} - 427049x^{18} - 812192x^{17} + 2322310x^{16} + 2957083x^{15} - 8374150x^{14} - 7010721x^{13} + 19892990x^{12} + 10323950x^{11} - 30317349x^{10} - 8509315x^9 + 28202264x^8 + 2963409x^7- 14792310x^6 + 220229x^5 + 3889959x^4 - 397376x^3 - 380960x^2 + 72460x - 2647, 2],[x^{27} - 3x^{26} - 53x^{25} + 160x^{24} + 1210x^{23} - 3680x^{22} - 15631x^{21} + 47937x^{20} + 126405x^{19} - 390929x^{18} - 670268x^{17} + 2085994x^{16} +       2381831x^{15} - 7410738x^{14} - 5717857x^{13} + 17541450x^{12} + 9224648x^{11} - 27292133x^{10} - 9719725x^9 + 27067090x^8 + 6187223x^7 - 16180684x^6 - 1937401x^5 + 5264635x^4 + 99800x^3 - 738836x^2 + 47264x + 14987, 2]$

\bibliography{bibgrau13}
\bibliographystyle{plain}

\end{document}